\documentclass[10pt]{article}
\usepackage{amsfonts}
\usepackage{natbib}
\usepackage{array}
\usepackage{amssymb}
\usepackage{amsmath}
\usepackage{graphicx,color}
\usepackage{amsfonts}
\usepackage[ignoreall]{geometry}

\newtheorem{theorem}{Theorem}

\newtheorem{lemma}[theorem]{Lemma}

\newtheorem{proposition}[theorem]{Proposition}
\newtheorem{remark}[theorem]{Remark}

\newenvironment{proof}[1][Proof]{\noindent\textbf{#1.} }{\ \rule{0.5em}{0.5em}}

\newcommand{\floor}[1] {\left\lfloor{#1}\right\rfloor}

\begin{document}


\renewcommand{\baselinestretch}{1.2}

\markright{
}

\markboth{\hfill{\footnotesize\rm FIRSTNAME1 LASTNAME1 AND FIRSTNAME2 LASTNAME2} \hfill}
{\hfill {\footnotesize\rm Sequential Analysis under Response Dependent Allocation} \hfill}

\renewcommand{\thefootnote}{}
$\ $\par


\fontsize{10.95}{14pt plus.8pt minus .6pt}\selectfont
\vspace{0.8pc}
\centerline{\large\bf Sequential Analysis of Cox  Model under}
\vspace{2pt}
\centerline{\large\bf  Response Dependent Allocation}
\vspace{.4cm}
\centerline{Xiaolong Luo, Gongjun Xu and Zhiliang Ying}
\vspace{.4cm}
\centerline{\it Celgene Corporation,
Columbia University and Columbia University}
\vspace{.55cm}
\fontsize{9}{11.5pt plus.8pt minus .6pt}\selectfont


\begin{quotation}
\noindent {\it Abstract:}
\cite{SellkeSiegmund} developed the Brownian approximation to the 
Cox partial likelihood score as a process of calendar time, laying the 
foundation for group sequential analysis of survival studies. We extend
their results to cover situations in which treatment allocations may depend
on observed outcomes. The new development makes use of the entry time and calendar time
along with the corresponding $\sigma$-filtrations to handle the natural information
accumulation. Large sample properties are established under suitable regularity conditions.

\vspace{9pt}
\par
\end{quotation}\par


\fontsize{10.95}{14pt plus.8pt minus .6pt}\selectfont

\section{Introduction}
The \cite{Cox}  proportional hazards model along with the partial likelihood \citep{Cox1975} has been extensively applied to survival data. 
The theoretical properties of the maximum partial likelihood estimator can be easily derived by expressing the partial likelihood score as a counting process based martingale integral; see \cite*{AndersonGill}, \cite{FlemingHarrington}, and \cite{KalbfleischPrentice}.

For sequential analysis, the partial likelihood score needs to be evaluated along the calendar time and its asymptotic behavior is crucial to deriving the corresponding group sequential methods.  Due to the staggered entry of patients, the partial likelihood score  as a process of calendar time is no longer a martingale integral. 
In a pioneering paper, \cite{SellkeSiegmund}  showed that the score process can still be approximated by the Brownian motion process, thereby laying the foundation for group sequential analysis of survival studies.  \cite{Slud} also established the Brownian approximation to the log-rank process for survival outcome under staggered entry. A Gaussian random
field approximation to the two-dimensional score process in the case
of two-sample comparison was established by  \cite{GuLai}; see also Andersen et al. (1993, Chapter 10). More general results about Gaussian random field approximation to the two-dimensional score process under the Cox proportional hazards
regression can be found in \cite*{BGY}, where modern empirical process theory is applied to derive certain key results, bypassing the martingale formulation.

The results of \cite{SellkeSiegmund} can be readily applied in the context of group sequential analysis as described in \cite{Pocock}, \cite{O'BrienFleming}, and \cite{LanDeMets}. However, their results are not applicable under adaptive designs where treatment allocation may depend on preceding outcomes. This is because the outcome variables are dependent so that neither the counting process-martingale argument nor the empirical process theory may be used to derive the desirable Brownian motion approximation. For some initial ideas of adaptive design, see \cite{Thompson} and \cite{Robbins}; for early works, see \cite{Zelen},  \cite{WeiDurham}, and \cite{Wei};  for more recent developments, see 
 \cite{FR1995} and \cite{HuRosenberger}.

The existing literature on response adaptive treatment allocation methods primarily deals with continuous or binary outcome variable. Recently  \cite{ZR2007} developed a parametric approach to survival outcomes. They assumed that survival times follow the exponential or, more generally, the Weibull family of distributions.
They showed that their approach can result in approximately optimal treatment allocation assuming survival times are relatively shorter than follow up period.

The main focus of this paper is to extend the results of \cite{SellkeSiegmund}  to the situation in which treatment allocations may depend on preceding outcomes. 
A key step in the new development is the expression of the
partial likelihood score process in terms of integrals over the
calendar and entry times.  As a result, the usual
martingale structure is preserved and can be applied to establish large sample properties. Indeed, it is shown that the partial likelihood score process is approximated by a time-rescaled Brownian motion process and that the maximum partial likelihood estimator is asymptotically normal.

The remainder of this paper is organized as follows. Section 2 first explains why the current martingale approach fails under the outcome dependent allocations, and then introduces a new approach. The corresponding functional central limit theorems are presented in Section 3, where convergence properties for the
corresponding maximum partial likelihood estimator are also
established. Some discussions are given in
Section 4. Most technical developments are presented in Appendix.

\section{Notation and model specification}

 We first introduce the setup and define some basic quantities. We will
consider a follow up study with calendar time period $[0, \tau]$, where
$\tau<\infty$. Let $n$ be the sample  size of the study. 
Denote by $U_{n,i}$ the entry time for individual $i$,
$i\geq 1$. For technical convenience, we assume throughout this paper
that the $U_{n,i}$ have no ties. Thus, without loss of generality, we
assume $U_{n,1}<U_{n,2}<\cdots<U_{n,i}<\cdots$. Define the associated counting
process for entry times
\begin{eqnarray}\label{enroll}
R_n(t)=\sum_{i\ge 1}I_{(U_{n,i}\leq t)}.
\end{eqnarray}
Note that $R_n(t)$ is the total number of enrollment up to time $t$ and $R_n(\tau)=n$. By large sample, we mean that $n$ goes to infinity while $\tau$ remains fixed. In other words, the situation considered here is high rate of entry over a fixed time period. 
An example of such kind in survival studies is the Beta-Blocker Heart Attack Trial \citep{BHAT}, where 3837 persons entered during the 27-month follow up period.
For notional convenience we shall henceforth omit subscript $n$ in $U_{n,i}$ whenever no confusion arises.

For subject $i$, let $T_i$ denote the survival time (since entry) and
$C_i$ the censoring time. Throughout the sequel, $a \wedge b = \min\{a, b\}$,
$a \vee b = \max\{a, b\}$, $a^+ = \max\{0, a\}$ and $a^- = \max\{0, -a\}$.
Let $\tilde{T}_i=T_i\wedge C_i$ and
$\Delta_i=I_{(T_i\leq C_i)}$, indicating failure (1) or censoring (0).
Thus, if $\Delta_i=1(0)$, then individual $i$ experiences failure
(censoring) at calendar time $U_i+\tilde{T}_i$. Furthermore, there
is a $p$-dimensional covariate vector
$Z_i$, which may include $i$th individual's
treatment assignment and certain relevant baseline characteristics.

\begin{figure}[h]
     \begin{center}
     \includegraphics[width=4in]{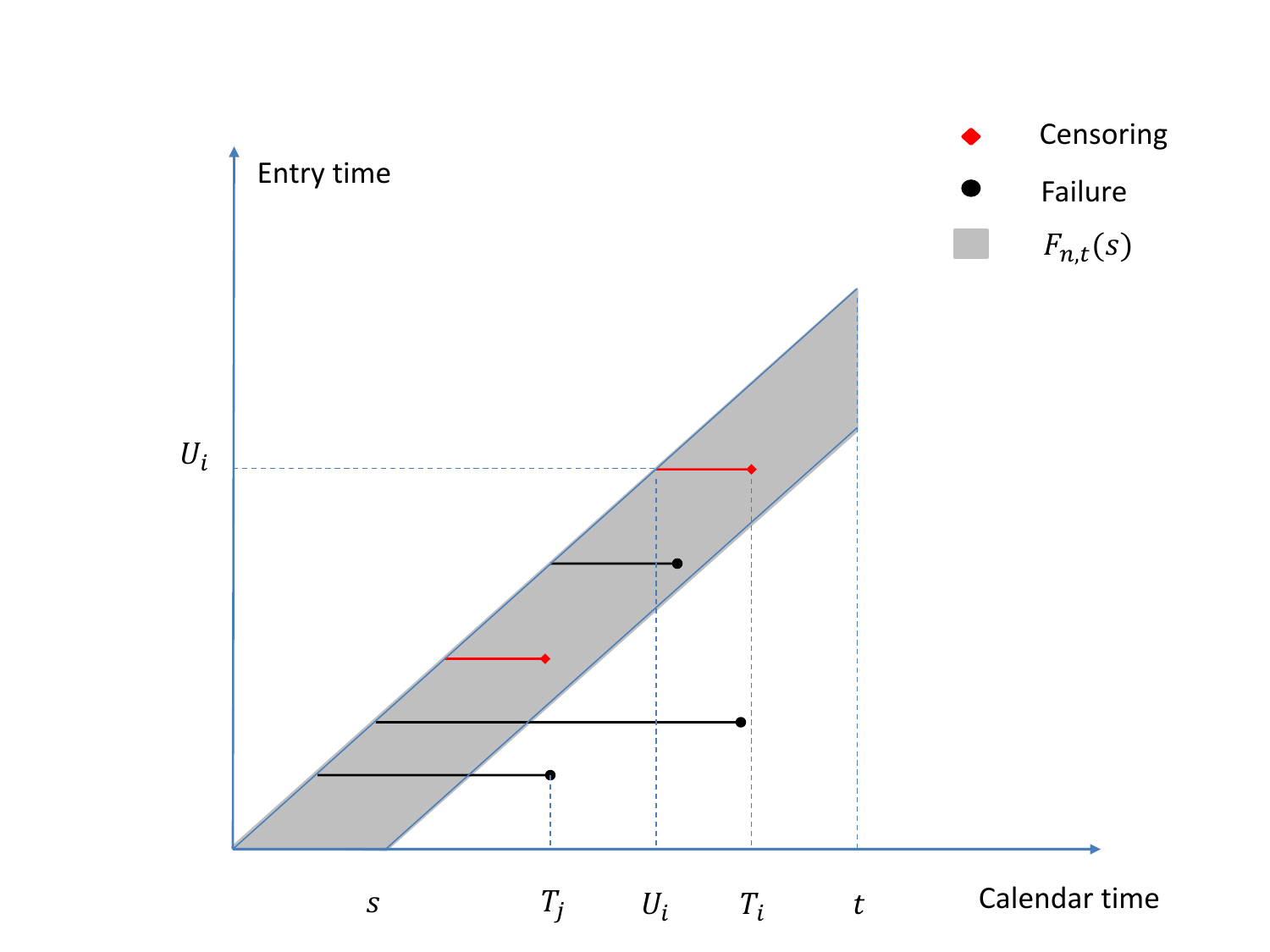}
      \caption{${\cal F}_{n,t}(s)$: $\sigma$-filtrations defined along calendar time and survival time. }
      \label{survival}
     \end{center}
\end{figure}

 We describe the Cox model specification with independent censoring under outcome dependent allocation as follow. For the $i$th subject, given $Z_i$, $T_i$ is conditionally independent of $C_i$ and $\{T_j$, $C_j$, $Z_j$; $j<i\}$ and has the following proportional hazards model specification
 $$\lambda_i(t)=\exp(\beta' Z_i)\lambda_0(t),$$
where $\beta$ is an unknown $p$-dimensional regression parameter of interest and $\lambda_0$ is the baseline hazard function. Note that under adaptive allocation, given $Z_j$, $T_j$ may not be independent of $T_{i}$ if $i>j$. This is because $Z_i$, which includes the treatment allocation of the $i$th subject, may depend on survival experiences of other subjects who enrolled before time $U_i$. For instance, in Figure \ref{survival}, we can see that $Z_i$ (and $T_i$) may depend on the survival information $T_j$ under the outcome dependent allocation scheme.  Compared with the independent enrollment scheme as in \cite{SellkeSiegmund}, where $\{T_i,C_i,Z_i\}$ are all assumed to be  independent,  outcome dependent allocation violates the independent assumption, raising the issue of validity for the existing sequential testing procedures. We will demonstrate the theoretical challenges arising from the violation of independence  in the next subsection, and propose our new approach in Subsection 2.2.

\subsection{Partial likelihood score process over survival time}

Under the usual nonadaptive allocation, i.e., observations from individual units are mutually independent, the partial likelihood \citep{Cox1975} takes form
\begin{equation}\label{PL}
PL(t) = \prod_{i:\, \tilde T_i\leq (t-U_i)^+,\atop \Delta_i=1}\biggr\{\left.{\exp(\beta'Z_i)}~\right/{\sum_{j:  \tilde T_j\leq (t-U_j)^+,\atop \tilde T_j\geq T_i}\exp(\beta'Z_j)}\biggr\}.
\end{equation}
Taking logarithm and differentiating with respect to $\beta$ result in the corresponding partial likelihood score process
\begin{eqnarray}\label{coxform}
 U_n(t)&=&\sum_{i:\,U_i\leq t}\int_0^t [Z_i - \bar Z_n (\beta; t,s)]N_i(t,ds),
\end{eqnarray}
where\begin{eqnarray}
{\bar Z}_n (\beta; t,s)&=&\frac{\sum_{i:\, U_i\leq t-s}^{}Z_i\exp(\beta' Z_i) I_{(\tilde T_i\geq s)}}
               {\sum_{i:\,U_i\leq t-s}\exp(\beta' Z_i)I_{(\tilde T_i\geq s)}}\nonumber
\end{eqnarray}
and $$N_i(t,s)=\Delta_iI_{(\tilde{T}_i\leq s\wedge (t-U_i)^+)}.$$
Let
\begin{eqnarray}\label{trandmrt}&&
 M_i(t,s)=N_i(t,s)
-\int_0^sI_{(\tilde{T}_i\wedge(t-U_i)^+\geq w)}\exp(\beta' Z_i)\lambda_0(w)dw.
\end{eqnarray}
It is well known that the partial likelihood score does not change numerically when the $N_i$ are replaced by the $M_i$, i.e.,
\begin{equation}\label{partiallik}
U_n(t)=\sum_{i:\,U_i\leq t}\int_0^t [Z_i - \bar Z_n (\beta; t,s)] M_i(t,ds).
\end{equation}

The integration in \eqref{partiallik} is with respect to survival time $s$. Under the usual independent sampling scheme, the $M_i$ are martingales as processes of $s$ with a suitably defined $\sigma$-filtration as in equation \eqref{Fs} below \citep{ABGK}. Furthermore, the integrands are predictable, so that $U_n(t)$ is a martingale integral with respect to survival time $s$. As a result, the martingale central limit theorem \citep{Rebolledo} can be applied to obtain the normal (Brownian) approximation.

Under the outcome dependent allocation, we now show that the martingale (along survival time $s$) argument is no longer valid.
For $s\geq 0$, let ${\cal F}_{n,t}(s)$ be the $\sigma$-filtration generated by observations up to survival time $s$ and calendar time $t$, i.e.,
\begin{eqnarray}\label{Fs}
{\cal F}_{n,t}(s)
&=&
\sigma\Big\{ I_{(U_i\leq t)},\quad  U_iI_{(U_i\leq t)},\quad Z_i I_{(U_i\leq t)},\\
&&~~~~ I_{(\tilde T_i\leq s\wedge (t-U_i)^+)},\quad
 N_i(t,s),\quad \tilde T_i I_{(\tilde T_i\leq s\wedge (t-U_i)^+)}
; \quad i=1,\cdots, n\Big\}.\nonumber
\end{eqnarray}
Figure \ref{survival} illustrates the information accumulated along survival time.
The grey trapezoid area shows the filtration ${\cal F}_{n,t}(s)$. From Figure \ref{survival}, we can see that for the $i$th subject enrolled at time $U_i$, although its survival time is less than $s$, its treatment allocation ($Z_i$) depends on the outcome information of $T_j$, which is outside of ${\cal F}_{n,t}(s)$.
Therefore, $ M_i(t,s)$ may not be a martingale with respect to filtration  ${\cal F}_{n,t}(s)$ under outcome dependent allocation.
However, if $\{T_i,C_i,Z_i\}$ are all independent as is the case in \cite{SellkeSiegmund} and \cite{GuLai}, the $ M_i(t,s)$ are still ${\cal F}_{n,t}(s)$ martingales in $s$ for any fixed $t$.

\subsection{Calendar time based score process}
In this subsection, we introduce a new way to represent the partial likelihood score so that a useful martingale structure will arise. The new representation expresses the score process in terms of integrals over entry time and calendar time. Use of entry time instead of survival time is natural in terms of the information accumulation from data and the adaptive treatment allocation process.

\begin{figure}
     \begin{center}
     \includegraphics[width=4in]{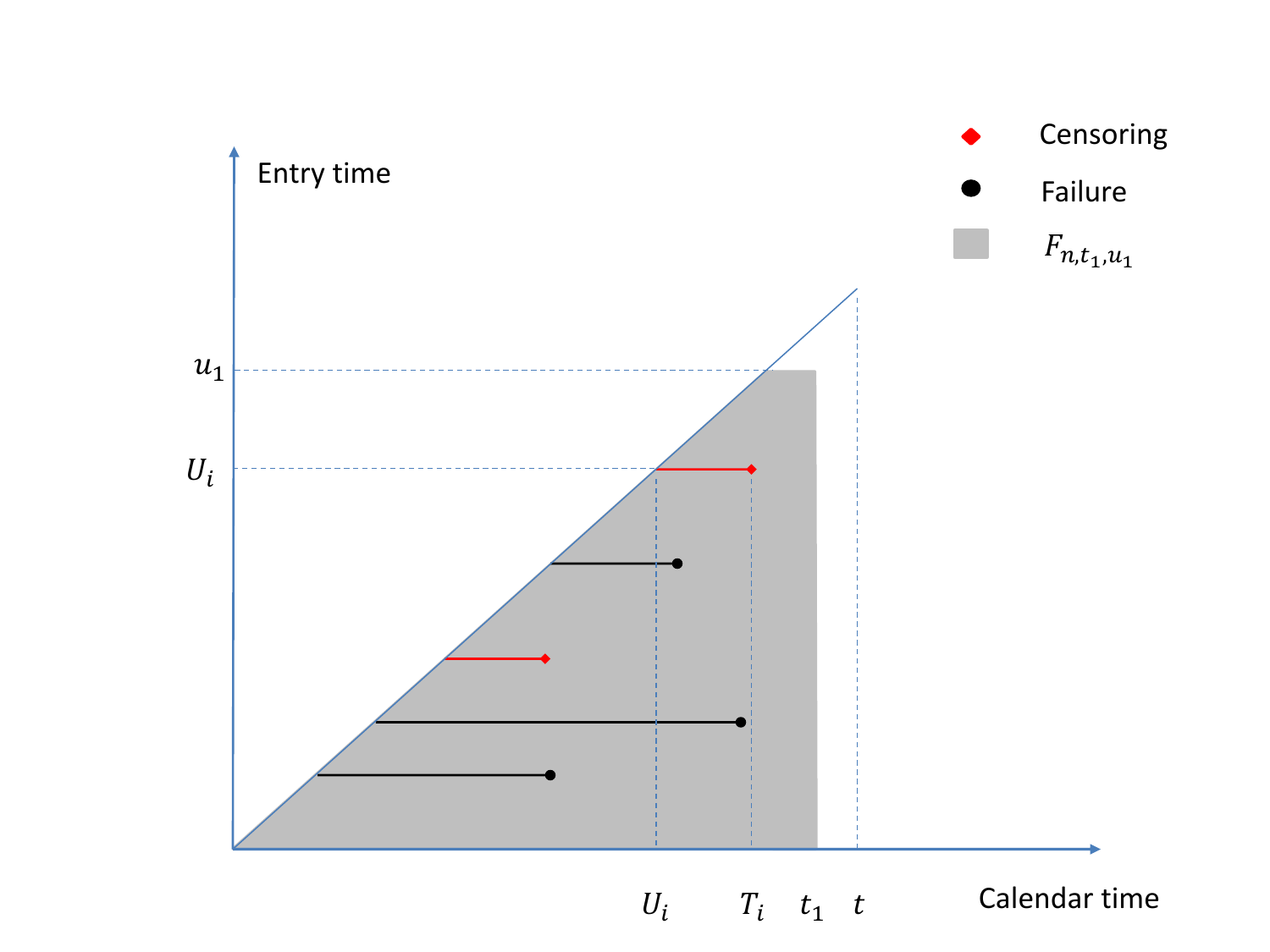}
     \caption{ ${\cal F}_{n,t,u}$: $\sigma$-filtrations defined along calendar time and entry time.}
     \label{F2}
     \end{center}
\end{figure}

With a slight abuse of notation, let $\tilde{T}_u$, $Z_u$, and
$\Delta_u$ refer to $\tilde{T}_i$, $Z_i$, and $\Delta_i$ when
$u=U_i$, which is well defined since the $U_i$ are distinct for
different $i$.
Define a random counting measure
\begin{eqnarray*}
p_n(ds\,du)
&=& I_{(u+\tilde{T}_u = s, \Delta_u=1)}dR(u),
\end{eqnarray*}
which defines a bivariate
counting process along both calendar time $s$ and entry time $u$. It equals 1 if there exists a subject $i$ such that $U_i=u$ and $T_i = s-u$; otherwise it equals 0. Based on the above two dimensional counting process,
the Cox score in \eqref{coxform} can be rewritten as an integral with respect to both calender time and entry time:
\begin{eqnarray}\label{calendar}
 U_n(t)&=&\int_0^t\int_0^s [Z_u - \bar Z_n(\beta; t,s-u)]p_n(ds\,du),
\end{eqnarray}
Let  ${\cal F}_{n,t}$  denote the corresponding $\sigma$-filtration containing all the information accumulation over calendar time period $[0,~ t]$, i.e.,
\begin{eqnarray*}
{\cal F}_{n,t}
&=& \sigma\Big\{ I_{(U_i\leq t)}, \quad
U_iI_{(U_i\leq t)},\quad Z_i I_{(U_i\leq t)},
\\
 &&~~~~ I_{(\tilde T_i\leq (t-U_i)^+)},\quad \Delta_i I_{(\tilde T_i\leq (t-U_i)^+)},\quad
 \tilde T_i I_{(\tilde T_i\leq  (t-U_i)^+)}; \quad i=1,\cdots, n\Big\}.\nonumber
\end{eqnarray*}
 A sub-$\sigma$-algebra of ${\cal F}_{n,t}$ that is of interest is defined by
\begin{eqnarray*}
{\cal F}_{n, t,\vartheta}
&=& \sigma\Big\{ I_{(U_i\leq \vartheta)},\quad U_iI_{(U_i\leq \vartheta)},\quad
Z_i I_{(U_i\leq \vartheta)},\quad I_{(\tilde T_i\leq (t-U_i)^+, U_i\leq \vartheta)},\\
 &&~~~~
\quad \Delta_i I_{(\tilde T_i\leq (t-U_i)^+,U_i\leq \vartheta)},\quad
\tilde T_i I_{(\tilde T_i\leq (t-U_i)^+,U_i\leq \vartheta)}
; \quad i=1,\cdots, n\Big\}.\nonumber
\end{eqnarray*}
Intuitively, ${\cal F}_{n, t,\vartheta}$ represents information up to calendar time $t$ for individuals who enrolled before time
$\vartheta$, where $0<\vartheta\leq t$. See Figure \ref{F2} for an illustration. The grey trapezoid area shows the filtration ${\cal F}_{n,t_1,u_1}$, which contains all the information up to calendar time $t_1$ and enrollment time $u_1$. Compared with Figure \ref{survival}, we can see that the treatment allocation information of the $i$th subject (enrolled at time $U_i$) is now included in the new filtration.

Without loss of generality, we shall assume throughout that $R(t)$
and $Z_t$ are predictable with respect to $\{{\cal F}_{n,t}, t\geq 0\}$, which is standard in survival analysis.
Note that for the $i$th subject, by the Dood-Meyer decomposition and the Cox model assumption, the compensator for the counting measure
$p(ds,u=U_i)= I_{(U_i+\tilde{T}_i = s, \Delta_u=1)}$ is $I_{(\Tilde{T}_i \ge s-U_i)}\exp\{\beta' Z_i\}\lambda_{0}(s-U_i)ds$ when $s>U_i$. This follows from the fact that 
$$P(U_i+T_i\in(s,s+ds] | {\cal F}_{n,s}) =  I_{(\Tilde{T}_i > s-U_i)}\exp\{\beta' Z_i\}\lambda_{0}(s-U_i)ds.$$
 More generally,
let
$$q_n(ds\,du)=I_{(\Tilde{T}_u \ge s-u)}\exp\{\beta' Z_u\}\lambda_{0}(s-u)dR(u)ds.$$
Note that $I_{(u<s)}q_n(ds\,du)$ is the compensator of $I_{(u<s)}p_n(ds\,du)$. Thus we have the following lemma.
\begin{lemma}\label{lemma1} For $t\in (0,\tau]$,
\begin{eqnarray}\label{martingale}
& & M_{n}(t)\triangleq \int_0^t\int_0^t I_{(u<s)}[p_n(ds\,du)-q_n(ds\,du)]
\end{eqnarray}
is a $\{{\cal F}_{n,t}, t\geq 0\}$ martingale. Moreover,
 for fixed $t$,
\begin{eqnarray}\label{martingale2way}
& & M_{n}(t,\vartheta)\triangleq
\int_0^t\int_0^{\vartheta}I_{(u<s)}[p_n(ds\,du)-q_n(ds\,du)],
\end{eqnarray}
as a process in $\vartheta$, is a $\{{\cal F}_{n,t,\vartheta},0\leq
\vartheta\leq t\}$ martingale.
\end{lemma}

Let  $$M_{n}(ds\,du)= I_{(u<s)}[p_n(ds\,du)-q_n(ds\,du)]$$ be the corresponding martingale measure.
The Cox score process in \eqref{calendar}  can then be written as
\begin{eqnarray*}
 U_n(\beta; t)&=&\int_0^t\int_0^t[ Z_u - \bar Z_n (\beta; t,s-u)]I_{(u<s)}p_n(ds\,du)\nonumber\\
 \label{coxscore}&=&\int_0^t\int_0^t[ Z_u - \bar Z_n (\beta; t,s-u)]M_{n}(ds\,du).
\end{eqnarray*}
More generally, we can define a two-parameter score process with
respect to calendar time $t$ and entry time $\vartheta$ as
\begin{eqnarray}
 U_n(\beta; t,\vartheta)=\int_0^t\int_0^{\vartheta}[ Z_u - \bar Z_n (\beta; t,s-u)]M_{n}(ds\,du).
\end{eqnarray}
Note that $U_n(\beta; t,t)=U_n(\beta; t)$.

The expression here for $U_n(\beta; t)$ is an
integral along the calendar time instead of the survival time as in standard counting process approach to survival analysis.  Through this framework, responses and
covariates history is expressed by the filtration ${\cal F}_{n,t}$. As a
result, it is not difficult to show that $M_{n,t}$ is a martingale with respect to $\sigma$-filtration ${\cal F}_{n,t}$ (Lemma \ref{lemma1}). This forms a crucial step for us to use the martingale central limit theorem to
obtain the convergence for $U_n(\beta; t)$; see Section
3 for more details.

\section{Large sample theory}

In this section, we establish large sample properties which are important for the usual statistical inferences, especially for sequential analysis. It is divided into two parts, with the first dealing with the score process and the second dealing with the estimator. 

\subsection{Weak convergence of score process}
The main effort of this subsection is to show the weak convergence of $U_n$ to a Gaussian random process. The result extends those of \cite{SellkeSiegmund}, \cite{GuLai}, and \cite{BGY} to cover the case with outcome dependent allocation schemes.

We adopt the setting of \cite{BGY} and restrict $t$ to $[0,\tau]$ with $\tau$ satisfying
\begin{eqnarray}\label{tail}
\liminf_{n\rightarrow\infty}\frac{1}{n}\sum_{i=1}^{n}P(\tilde T_i\geq \tau)>0
\end{eqnarray}
and $\lambda_0$ being bounded on $[0,\tau]$. This entails that here we are adopting the asymptotics in terms of a high rate of entry over a fixed time interval (large $n$) as opposed to a fixed rate of entry over a long time interval (large $\tau$); see Siegmund (1985; p. 126). It allows us to develop a Gaussian random field approximation as in \cite{GuLai}, which also assumes large $n$. For asymptotics under large $\tau$, certain rescaling is needed and the corresponding Gaussian approximations may also be developed under certain stability assumptions \citep{Siegmund1985}.





For a $p$-dimensional covariate vector $Z$ with regression parameter vector
$\beta$, let $Z^{\otimes 0}=1$, $Z^{\otimes
1}=Z$, and $Z^{\otimes 2}=ZZ'$.
For $k=0,1$ and $2$,  $\vartheta>0$ and $w>0$, let
\begin{eqnarray}\label{gamma}
S_{n,k}(\beta; \vartheta,w) &=&
 \sum_{i:\,U_i\leq \vartheta}Z_i^{\otimes k}\exp(\beta' Z_i)I_{(\tilde{T}_i\geq w)}\nonumber\\
&=&
 \int_0^{\vartheta} I_{(\Tilde{T}_u\ge w)} Z^{\otimes k}_u\exp(\beta' Z_u)dR(u).
\end{eqnarray}
Recall the score processes defined as in Section 2.2:
\begin{eqnarray}
 U_n(\beta; t)&=&\int_0^t\int_0^t[ Z_u - \bar Z_n (\beta; t,s-u)]M_{n}(ds\,du),\nonumber\\
 U_n(\beta; t,\vartheta)&=&\int_0^t\int_0^{\vartheta}[Z_u - {\bar Z}_n (\beta;t,s-u)]M_{n}(ds\,du),\nonumber
\end{eqnarray}
where
\begin{eqnarray}
{\bar Z}_n (\beta; t,w)&=&\frac{S_{n,1}(\beta; t-w,w) }
               {S_{n,0}(\beta; t-w,w)}.\nonumber
\end{eqnarray}

Let $\beta_0$ be  the true regression parameter.
The following conditions are needed for our main results. 
\begin{itemize}
\item[C1] The $Z_i$ are uniformly bounded in the sense that
there exists a non-random constant $K_\tau$ such that,
$\sup_{i:\,U_i\leq \tau} |Z_i|\leq K_{\tau},$ where $|\cdot|$ denotes the $L_1$ norm for a $p$-dimensional vector.

\item[C2]  For $k=0$, $1$ and $2$, there exist non-random constants
$\bar E_k(\vartheta,w)$ such that as $n\rightarrow
\infty$,
$$\frac{1}{n}S_{n,k}(\beta_0; \vartheta,w)- \bar E_k(\vartheta,w)\xrightarrow[]{P}
0,$$ uniformly for all positive $\vartheta,w $ satisfying $\vartheta+w\leq
\tau$.

\end{itemize}

\begin{remark}
Conditions C1 and C2 are analogous to Conditions 1-3 in \cite{BGY}. In particular, C1 may be extended to a moment condition on $Z$ for the components related to the baseline covariates. Condition C2 is required so that the sample moments for the $Z_i$ are stable.
\end{remark}


We now state the main result of this subsection.
\begin{theorem}\label{mle}
Suppose that Conditions $C1$ and $C2$ are satisfied. Then the following convergence results hold for the partial likelihood score processes $U_n(\beta_0;t)$ and $U_n(\beta_0;t,\vartheta)$. \\
$(i)$ $n^{-1/2}U_n(\beta_0;t)$ converges
weakly to a vector-valued zero-mean Gaussian process $\xi$ on
$[0,\tau]$ with covariance function
\begin{equation*}
E[\xi(t_1)\xi'(t_2)]=\int_0^{t_1\wedge t_2}\left[\bar E_2(t_1\wedge t_2-w,w)
-\frac{\bar E_1^{\otimes 2}(t_1\wedge t_2-w,w)}
{\bar E_0(t_1\wedge t_2-w,w)}\right]\lambda_0(w) dw.
\end{equation*}
$(ii)$ $n^{-1/2}U_n(\beta_0;t,\vartheta)$ converges
weakly to a vector-valued zero-mean Gaussian random field $\tilde
\xi(t,\vartheta)$ on $\{(t,\vartheta): 0\leq \vartheta\leq t\leq
\tau\}$ with  covariance function
\begin{eqnarray*}
&&E[\tilde \xi(t_1,u_1)\tilde\xi'(t_2,u_2)]\\
&=&\int_0^{t_1\wedge t_2}\biggr[\tilde E_2(u_1, u_2, t_1, t_2, w) - 
2\frac{\bar E_1(t_1\wedge t_2-w,w)}{\bar E_0(t_1\wedge t_2-w,w)}{\tilde E_1'(u_1,u_2,t_1,t_2,w)}\\
&&
~~~~~~+\frac{\bar E_1^{\otimes 2}(t_1\wedge t_2-w,w)}{(\bar E_0(t_1\wedge t_2-w,w))^2}{\tilde E_0(u_1,u_2,t_1,t_2,w)}\biggr]\lambda_0(w)dw,
\end{eqnarray*}
where $\tilde E_k(u_1, u_2, t_1, t_2, w) = \bar E_k(u_1\wedge u_2\wedge(t_1\wedge t_2-w),w)$, $k=0, 1$ and $2$.
\end{theorem}

\begin{remark}
Theorem \ref{mle} extends the existing results by allowing allocation schemes to be depent on previous information. In addition,  it implies that $\tilde\xi$ has independent increments in  calender time $t$. Thus the diagonal process $\tilde\xi(t,t)=\xi(t)$ is a time-rescaled Brownian motion when $dim(Z) = 1$, and a vector-valued Gaussian process with independent increments when $dim(Z) > 1$.
\end{remark}

\begin{remark}
To apply Theorem \ref{mle}, we need to estimate the covariance function $E[\tilde \xi(t_1,u_1)\tilde\xi'(t_2,u_2)]$. A natural approach is to replace the unknown quantities $\bar E_{k}$ and $\Lambda(\cdot)$ with $S_{n,k}/n$ and the Nelson-Aalen estimator respectively. Consistency of the corresponding covariance estimator can be  derived under Conditions C1 and C2. 
\end{remark}

The following lemma plays a key role in the proof of Theorem 4. Its proof is given in the Appendix.

\begin{lemma}\label{lemma4'}
Under the same assumptions as those of
Theorem \ref{mle}, we have
\begin{eqnarray}
\sup_{\vartheta,t\in[0,\tau]}\frac{1}{{\sqrt{n}}}\int_0^{t}\int_0^{\vartheta}\left[\bar Z_n (\beta_0; t,s-u)-\frac{\bar E_1(t-(s-u),s-u)}{\bar E_0(t-(s-u),s-u)}\right]M_{n}(ds\,du)
\xrightarrow[]{P} 0.\nonumber
\end{eqnarray}
\end{lemma}

\begin{remark}
 Lemma \ref{lemma4'} shows that $\bar Z_n$ may be replaced by its (non-random) limit. The replacement makes it easy to use the martingale structure along the calendar time and the entry time without appealing to  the empirical process theory which may not be applicable under outcome dependent  allocation schemes.
\end{remark}

\begin{proof}[Proof of Theorem \ref{mle}]
Note that when $\vartheta=t$, $U_n(\beta_0;t,\vartheta)=U_n(\beta_0; t)$. So we only need to prove the weak convergence of $n^{-1/2}U_n(\beta_0;t,\vartheta)$.
 By Lemma \ref{lemma4'}, it suffices to show the weak convergence of
 $$n^{-1/2}\tilde U_n(\beta_0; t,\vartheta)=n^{-1/2}\int_0^t\int_0^{\vartheta}
 \left[Z_u - \frac{ \bar E_1(t-(s-u),s-u)}{ \bar E_0(t-(s-u),s-u)}\right]M_{n}(ds\,du).$$
We first show that for any positive integer
$k$ and partition $0\leq u_1<\cdots<u_{k}\leq \tau$,
$\{n^{-1/2}\tilde{U}_n(
\beta_0;t,u_1)$,$\cdots$,$n^{-1/2}\tilde{U}_n(
\beta_0;t,u_{k})$, $0\leq t\leq \tau\}$ converges weakly to a
multivariate Gaussian process $\{\tilde \xi(t,u_1)$,$\cdots$,
$\tilde \xi(t,u_{k}), 0\leq t\leq \tau\}$.
By Lemma \ref{lemma1},
we have that
$\{\tilde{U}_n(\beta_0;t,u_j),{\cal F}_{n,t}, 0\leq t\leq
\tau\}$ are martingales along calendar time $t$ with predictable variation processes
 \begin{eqnarray*}
&&\langle n^{-1/2}\tilde U_n(\beta_0; \cdot,u_i),n^{-1/2}\tilde U_n(\beta_0; \cdot,u_j)\rangle(t)\nonumber\\
&=&
\frac{1}{n}\int_0^t\int_0^{u_i\wedge u_j\wedge s}
\left[Z_u-\frac{ \bar E_1(t-(s-u),s-u)}{ \bar E_0(t-(s-u),s-u)}\right]^{\otimes2}
q_n(ds\,du)\nonumber\\
&\xrightarrow[]{P}&
\int_0^{t}\biggr[\bar E_2(u_i\wedge u_j\wedge(t-w), w) -
2\frac{\bar E_1(t-w,w)}{\bar E_0(t-w,w)}{\bar E_1'(u_i\wedge u_j\wedge(t-w),w)}\\
&&
~~~~~~+\frac{\bar E_1^{\otimes 2}(t-w,w)}{(\bar E_0(t-w,w))^2}{\bar E_0(u_i\wedge u_j\wedge(t-w),w)}\biggr]\lambda_0(w)dw\nonumber,
\end{eqnarray*}
where the convergence in probability is uniform in $t$ and follows
from Condition $C2$.
By the martingale central limit theorem \citep{Rebolledo}, we obtain that
 any linear combination of $\{n^{-1/2}\tilde{U}_n(
\beta_0;t,u_1)$, $\cdots$, $n^{-1/2}\tilde{U}_n(
\beta_0;t,u_{k})$, $0\leq t\leq \tau\}$ converges weakly to the corresponding linear transformation of
$\{\tilde \xi(t,u_1)$, $\cdots$,
$\tilde \xi(t,u_{k}), 0\leq t\leq \tau\}$. Therefore, we obtain
the weak convergence of $\{n^{-1/2}$ $\tilde{U}_n(
\beta_0;t,u_1)$, $\cdots$, $n^{-1/2}\tilde{U}_n$ $(
\beta_0;t,u_{k})$, $0\leq t\leq \tau\}$ 
via the Cram\'{e}r-Wold device. In particular, $n^{-1/2}\tilde U_n(\beta_0; t,\vartheta)$
converges in finite dimensional distribution to a Gaussian random field. 

In the Appendix, it will be shown (Proposition \ref{prop1}) that for any $\epsilon>0$, there exist a
constant $k_0<\infty$ and partition $0=u_{n,0}\leq
u_{n,1}\leq\cdots\leq u_{n,k_0}=\tau$ such
that for all large $n$,
$$P\left( \max_{0\leq j<k_0} \sup_{\vartheta\in [u_{n,j},u_{n,j+1}];\atop 0\leq t\leq \tau}\frac{1}{\sqrt{n}}|\tilde U_{n}(\beta_0; t,\vartheta) - \tilde U_{n}(\beta_0; t,u_{n,j})|\geq \epsilon\right)\leq \epsilon.$$
Thus, $n^{-1/2}\tilde U_n(\beta_0; t,\vartheta)$ is tight. Combing this with the above finite dimensional distributional convergence result, we obtain the desired conclusion.
\end{proof}

\subsection{Asymptotic normality of maximum partial likelihood estimator}
From Section 3.1, we know that $U_n(\beta_0;t,\vartheta)$ converges to a zero mean Gaussian process in
the direction of both $t$ and $\vartheta$. Thus, we may use it to obtain an asymptotically unbiased estimator of $\beta$ for each fixed $(t, \vartheta)$. Specifically, let 
$\hat{\beta}(t,\vartheta)$ be the
solution to $U_n(\beta;t,\vartheta)=0$. Note that at $\vartheta=t$, $\hat{\beta}(t,t)$ is simply the maximum partial likelihood estimator with observable data at calendar time $t$. We will show in this subsection that $\hat\beta(t,\vartheta)$ is asymptotically normal, or, more precisely, $\sqrt n (\hat\beta(t,\vartheta)-\beta_0)$ converges  weakly to a zero mean Gaussian process. 

To establish the asymptotic normality, we first state the following condition, which ensures that the information matrix is nonsingular when normalized by the sample size $n$.

\medskip
{C3}. There exists $\tau_0\in (0,\tau]$ such that for all $(\vartheta,\tau)$ satisfying $\tau_0\leq \vartheta\leq t\leq \tau$,
\begin{eqnarray*}
&&\lambda_{min} (A(t,\vartheta))\geq v_0>0, ~a.s.,
\end{eqnarray*}
where 
\begin{eqnarray*}
A(t,\vartheta)=\int_0^{t}\left[\bar E_2(t-w,w)
-\frac{\bar E_1^{\otimes 2}(t-w,w)}
{\bar E_0(t-w,w)}\right]\frac{\bar E_0(\vartheta\wedge(t-w),w)}{\bar E_0(t-w,w)}
\lambda_0(w) dw,
\end{eqnarray*}
$\bar E_{k}$ is defined as in Condition C2 and
$\lambda_{min}(A)$ denotes the minimum eigenvalue of a symmetric
matrix A.

\begin{theorem}\label{consistency}
 Suppose that Conditions $C1$, $C2$, and $C3$ are satisfied. Then, 
$\{\sqrt{n}(\hat{\beta}(t,\vartheta)-\beta_0),
\tau_0\leq \vartheta\leq t\leq \tau\}$ converges weakly to a
vector-valued zero-mean Gaussian process $\eta$ with covariance function
$$E[\eta(t_1,u_1)\eta'(t_2,u_2)]=\left(A(t_1,u_1)\right)^{-1}E[\tilde\xi(t_1,u_1)\tilde\xi'(t_2,u_2)]
\left(A(t_2,u_2)\right)^{-1},$$
where $\tilde\xi$ is the Gaussian process defined as in Theorem
\ref{mle}.
\end{theorem}

\begin{proof}[Proof of Theorem \ref{consistency}]
By Lenglart's inequality (Lemma \ref{lenglart}), we have that as
$n\rightarrow\infty$, 
\begin{equation}\label{lem9c1}
\sup_{0\leq \vartheta\leq t\leq \tau}
\Big|\frac{1}{n}U_n(\beta_0;t,\vartheta)\Big|\xrightarrow[]{P}0.
\end{equation}
Condition C2 implies that 
\begin{eqnarray}\label{4.4}
&& \sup_{0\leq\vartheta\leq t\leq\tau}\biggr|\frac{1}{n}\frac{\partial}{\partial
\beta} U_n(\beta_0;t,\vartheta) + A(t,\vartheta)\biggr|\xrightarrow[]{P}0.
\end{eqnarray}
Since $\frac{1}{n}\frac{\partial}{\partial \beta}
U(\beta_0;t,\vartheta)$ has a uniformly bounded
derivative with respect to $\beta$, Condition C3 and
\eqref{4.4} imply that there exists a neighborhood of
$\beta_0$, ${\cal N}(\beta_0)$, such that
\begin{eqnarray}\label{neigb}
\liminf_{n\rightarrow\infty}\inf_{\tau_0\leq \vartheta\leq t\leq \tau}
\inf_{\beta\in {\cal N}(\beta_0)}
\lambda_{min}\Big(-\frac{1}{n}\frac{\partial}{\partial\beta}
U(\beta;t,\vartheta)\Big)\geq\frac{v_0}{2}>0.
\end{eqnarray} 
Therefore, by \eqref{lem9c1}, \eqref{neigb}, and Lemma \ref{bgy} in the Appendix (A4), together with the positive definiteness of $-\frac{1}{n}\frac{\partial}{\partial \beta}
U(\beta;t,\vartheta)$,  we obtain the uniform consistence of $\hat\beta(t,\vartheta)$, that is
$$\sup_{\tau_0\leq \vartheta\leq t\leq\tau}|\hat\beta(t,\vartheta)-\beta_0|\xrightarrow[]{P} 0.$$

By the Taylor series expansion, we have that
\begin{eqnarray*}
0&=&\frac{1}{\sqrt{n}}U(\hat{\beta}(t,\vartheta);t,\vartheta)\\
&=&\frac{1}{\sqrt{n}}U(\beta_0;t,\vartheta)+\frac{1}{n}\frac{\partial}{\partial
\beta}
U(\beta_0;t,\vartheta)\sqrt{n}(\hat{\beta}(t,\vartheta)-\beta_0)+o_p(1),
\end{eqnarray*}
uniformly in $\tau_0\leq \vartheta\leq t\leq \tau$. Therefore, 
$$\sqrt{n}(\hat\beta(t,\vartheta)-\beta_0) = 
-\left(\frac{1}{n}\frac{\partial}{\partial
\beta}
U(\beta_0;t,\vartheta)\right)^{-1}\frac{1}{\sqrt{n}}U(\beta_0;t,\vartheta)
+o_p(1).$$
The weak convergence of
$\sqrt{n}(\hat{ \beta}(t,\vartheta)-\beta_0)$
follows from the above expansion and Theorem \ref{mle}. 
\end{proof}

\section{Discussion}
This paper considers the Cox model based sequential analysis of survival studies when treatment allocation may depend on survival outcomes observed prior to the time of treatment assignment. It develops an
approach based on calendar time and entry time filtrations to obtain basic martingales and to bypass independence assumption that may likely be violated. 
Desirable asymptotic properties are then obtained under suitable regularity conditions, showing that the Brownian approximation to the partial likelihood score process, initially developed by \cite{SellkeSiegmund},  is still valid. As a consequence,  the usual group sequential boundaries such as those by  \cite{Pocock}, \cite{O'BrienFleming}, and \cite{LanDeMets} can be used. As an example, suppose that treatment assignments in later stages of the Beta-Blocker Heart Attack Trial (BHAT, 1982) were adapted to the survival outcomes so that more patients are to be allocated to the treatment arm if that arm has significantly more favorable results. Then the Brownian approximation would still be valid for deriving the corresponding group sequential boundary.

One of the limitations of the asymptotic theory developed here is the assumption of high accrual rate in a fixed follow up period. Such an assumption entails that a significant portion of survival experiences from previously entered subjects may not be fully available for optimal treatment allocation due to delayed survival outcomes. Consequently, the asymptotically optimal treatment allocation ratio as discussed in  \cite{ZR2007}  may not be attainable. On the other hand, the flexibility of using all observed survival outcomes could alleviate this deficiency of delayed response. 

Another way to formulate large sample setting is to assume large time, rather than high accrual rate. That is to consider $\tau$ (follow up period) going to infinity. Under this formulation, for large (calendar time) $t$, the proportion of observed outcomes from previously entered subject will tend to $1$ as $t$ goes to infinity, making the asymptotically optimal treatment allocation feasible. When there is no other explanatory variable besides a dichotomous treatment allocation, it is not difficult to extend the present approach by rescaling of time through the ``compensator''.  In general, it may require additional assumptions on the explanatory variables in order to establish the vector-valued Gaussian martingale approximation to the multivariate score process. 

To alleviate the effect of delayed survival outcomes, certain surrogate variables (markers) for the survival time may be used for the purpose of treatment allocation. For example, in the BATTLE trial  \citep{Zhou2008,Kim2011}, if patients' survival times were the endpoint, then one could use progression-free survival as a surrogate variable. It will be of interest to develop a similar theoretic framework under which the Brownian approximation may be used.  

The approach developed here may be extended to other follow-up studies with more general
outcome variables. For studies with longitudinal outcomes, dynamic
regression models have been proposed and studied \citep{MartinussenScheike}. Adaptive and outcome dependent designs for such studies may result in staggered entry and dependent observation units. We believe the general approach developed in this paper can be extended to deal with such designs.

\section*{Appendix.}
\subsection*{A1: Proof of Lemma \ref{lemma4'}}
To simplify notation, we let
$f_{n}(t,s,u)=\bar Z_n (\beta_0; t,s-u)$ and 
$$g(t,s,u)=\frac{\bar E_1(t-(s-u),s-u)}{\bar E_0(t-(s-u),s-u)}.$$
For a subject who enrolled into the study at time $u$, define, for $s\in [u,\tau]$, counting measure
$$p_{n,u}(ds)= I(u+\Tilde{T}_u= s).$$
Under the $\sigma$-filtration ${\cal F}_{n,t}$, it is easy to see the compensator for $p_{n,u}(ds)$ is 
$$q_{n,u}(ds)= I_{(\tilde T_u\geq s-u)}\exp\{\beta'Z_u\}\lambda_0(s-u)ds.$$
Thus,
$$M_{n,u}(ds)= I_{(u<s)}[p_{n,u}(ds)-q_{n,u}(ds)]$$
is a martingale measure. Comparing this with (\ref{trandmrt}), it follows that $M_i(t,s)$ is a martingale as a process in $s$, since for the $i$th subject with entry time $u=U_i$, $M_{n,U_i}(ds)=I(U_i+\Tilde{T}_i\in
ds)-q_{n,U_i}(ds)=M_{i}(t,d(s-U_i))$ for $s>U_i$.
 Define martingale integral
$$M_{n,u}(t)=\int_u^t M_{n,u}(ds),$$
which is the total measure on interval $[u,t]$. Let
 $$M_{n}(t,du)=\left[\int_u^t M_{n,u}(ds)\right] dR(u)=M_{n,u}(t) dR(u),$$
which defines a random measure along entry time for subjects who enrolled into the study before time $t$.

 Under the above notation,  for $M_{n}(t,\vartheta)$ defined as in \eqref{martingale2way}, we have the following identity 
 $$M_{n}(t,\vartheta)=\int_0^t\int_0^{\vartheta} M_n(ds\,du)=\int_0^\vartheta
 M_{n}(t,du).$$
Note that from Lemma \ref{lemma1}, $M_{n}(t,\vartheta)$ is  a
martingale along both calendar and entry times, i.e.,
 $M_{n}(t,\vartheta)$ is a martingale in $t$ for any fixed $\vartheta$ and a martingale in $\vartheta$ for any $t$. When $\vartheta=t$, we have
  $M_{n}(t,t)=\int_0^t\int_0^t M_{n}(ds\,du)$, which is $M_{n}(t)$. Similarly,
define random integral $\tilde{M}_{n}(w,\vartheta)$ with respect to
survival time $w$ and entry time $\vartheta$ by
\begin{equation}\label{tildeM}
\tilde{M}_{n}(w,\vartheta)=
\int_0^{\vartheta} M_{n}(w+u, du)\left
(=\int_0^{\vartheta} M_{n,u}(w+u)dR(u)\right).
\end{equation}
Note that $\tilde{M}_{n}(w,\vartheta)$ is defined on the information
observed before entry time $\vartheta$ and survival time $w$.

To prove Lemma \ref{lemma4'}, we need the following two
propositions, whose proofs are given in
subsections A2 and A3, respectively. Proposition \ref{prop1}
shows that $M_{n}(t,\vartheta)/\sqrt{n}$ is tight along
calendar and entry times while
Proposition \ref{prop2} shows the tightness property for
$\tilde{M}_{n}(w,\vartheta)/\sqrt{n}$ along survival and entry
times.

\begin{proposition}\label{prop1}
Under Conditions C1 and C2, for any $\epsilon>0$, there exist 
constant $n_0<\infty$ and partition $0=u_{n,0}\leq
u_{n,1}\leq\cdots\leq u_{n,n_0}=\tau$ such
that for all large n,
$$P\left( \max_{0\leq j<n_0} \sup_{\vartheta\in [u_{n,j},u_{n,j+1}];\atop 0\leq t\leq \tau}|W_{n,t,\vartheta}-W_{n,t,u_{n,j}}|\geq \epsilon\right)\leq \epsilon,$$
where $W_{n,t,\vartheta}=M_{n}(t,\vartheta)/\sqrt{n}$.
\end{proposition}

\begin{proposition}\label{prop2}
Under Conditions C1 and C2, for any $\epsilon>0$, there exist
partitions $0=w_0<w_1<\cdots<w_{N_0}=\tau$ and $0=u_{n,0}\leq
u_{n,1}\leq \cdots\leq u_{n,n_0}=\tau$   such that for all large n,
$$P\left( \max_{0\leq j<n_0 \atop 0\leq k<N_0 } \sup_{\vartheta\in [u_{n,j},u_{n,j+1}]\atop w\in[w_k,w_{k+1}]}|\tilde{W}_{n,w,\vartheta}-\tilde{W}_{n,w_k,u_{n,j}}|\geq \epsilon\right)\leq \epsilon,$$
where
$\tilde{W}_{n,w,\vartheta}=\tilde{M}_{n}(w,\vartheta)/\sqrt{n}$.
\end{proposition}

\begin{proof}[Proof of Lemma \ref{lemma4'}]
For any $(\vartheta,t)$ such that $0\leq \vartheta\leq t\leq \tau$, by changing the integration order, we have that
\begin{eqnarray}\label{proof10}
&& \frac{1}{\sqrt{n} }\int_0^t\int_0^{\vartheta}(f_n(t,s,u)-g(t,s,u))M_n(ds\,du)  \nonumber \\
& = &\frac{1}{\sqrt{n} }\int_0^t\int_0^{s\wedge \vartheta}(f_n(t,s,u)-g(t,s,u))(p_n(ds\,du)-q_n(ds\,du))  \nonumber \\
& =&\frac{1}{\sqrt{n} }  \int_0^{\vartheta}\int_u^t(f_n(t,s,u)-g(t,s,u))(p_n(du\,ds)-q_n(du\,ds))   \nonumber \\
& =&\frac{1}{\sqrt{n} }
\int_0^\vartheta\left[\int_u^t(f_n(t,s,u)-g(t,s,u))M_{n,u}(ds)\right]dR(u)
\nonumber\\
 &=&
\frac{1}{\sqrt{n}}\int_0^\vartheta\biggr[M_{n,u}(t)(f_n(t,t,u)-g(t,t,u))
-M_{n,u}(u)(f_n(t,u,u)-g(t,u,u)) \nonumber \\
 &&-\int_u^t M_{n,u}(s)(f_n(t,ds,u)-g(t,ds,u))\biggr]dR(u) +o_p(1),  
\end{eqnarray}
where the last equation follows from the integration-by-parts formula. Inclusion of $o_p(1)$ is due to the discontinuity of both the integrand and the integrator functions when the integration-by-parts formula is used.
Therefore,  by the definition of $M_{n}(t,du)$ and the fact that $M_{n,u}(u)=0$, 
we get
\begin{eqnarray}\label{proof1}
  \eqref{proof10} 
&=&   \frac{1}{\sqrt{n}} \int_0^\vartheta(f_n(t,t,u)-g(t,t,u))M_{n}(t,du)\nonumber\\
&  & - \frac{1}{\sqrt{n}}\int_0^\vartheta\left[\int_u^t
M_{n,u}(s)(f_n(t,ds,u)-g(t,ds,u))\right]dR(u)+o_p(1).\notag\\
\end{eqnarray}
In view of \eqref{proof1}, it remains to show that the two leading terms in (\ref{proof1}) are negligible. 

For the first term, taking integration by parts, we have that
\begin{eqnarray}\label{14}
&  & \left|\frac{1}{\sqrt{n}}\int_0^\vartheta(f_n(t,t,u)-g(t,t,u))M_{n}(t,du)\right|   \nonumber \\
& =& \biggr|\frac{1}{\sqrt{n}}(f_n(t,t,\vartheta)-g(t,t,\vartheta))M_{n}(t,\vartheta)
-\frac{1}{\sqrt{n}}(f_n(t,t,0)-g(t,t,0))M_{n}(t,0) \nonumber\\
&&- \frac{1}{\sqrt{n}}\int_0^\vartheta M_{n}(t,u)(f_n(t,t,du)-g(t,t,du)) \biggr| +o_p(1)\nonumber\\
& \leq & \left|\frac{1}{\sqrt{n}}(f_n(t,t,\vartheta)-g(t,t,\vartheta))M_{n}(t,\vartheta) \right| \notag\\
&&+ \frac{1}{\sqrt{n}}\left|\int_0^\vartheta M_{n}(t,u)(f_n(t,t,du)-g(t,t,du))\right|+o_p(1).
\end{eqnarray}
From Proposition \ref{prop1}, we have, for any $\epsilon>0$, there exists a partition
$0=u_{n,0}\leq u_{n,1}\leq\cdots\leq u_{n,n_0}=\tau$ such that for all
large $n$, 
$$P\left (\sup_{i; u\in
(u_{n,i},u_{n,i+1}]}
\frac{1}{\sqrt{n}}|M_{n}(t,u_{n,i+1})-M_{n}(t,u)|<\epsilon\right)\ge 1-\epsilon.$$ 
Combining this with \eqref{14},  for all large $n$, the following result holds uniformly
on $0\leq \vartheta\leq t\leq \tau$ with
probability at least $1-2\epsilon$:
\begin{eqnarray}\label{20}
\eqref{14}
& \leq & \left|\frac{1}{\sqrt{n}}(f_n(t,t,\vartheta)-g(t,t,\vartheta))M_{n}(t,\vartheta) \right|\nonumber\\
&  & + \frac{1}{\sqrt{n}}\sum_{i=1}^{n_0}\left|\int_{u_{n,i-1}}^{u_{n,i}}M_{n}(t,u_{n,i})
(f_n(t,t,du)-g(t,t,du))\right| + 2\epsilon K\nonumber\\
& \leq & 3\epsilon K,
\end{eqnarray}
where $K$ is the total variation bound for
$f_n(t,s,u)(=f_n(t,s-u,0))$, and the last inequality follows from 
Lenglart's inequality (Lemma \ref{lenglart}). Since $\epsilon$ can be arbitrarily small, the first term is negligible.  

For the second term, by the definitions of $f_n$ and $g$, we have
$f_n(t,s,u)=f_n(t,s-u,0)$ and $g(t,s,u)=g(t,s-u,0)$. Therefore, 
\begin{eqnarray}\label{proof2}
&  &\frac{1}{\sqrt{n}} \int_0^\vartheta\int_u^t M_{n,u}(s)(f_n(t,ds,u)-g(t,ds,u))dR(u) \nonumber\\
& = &\frac{1}{\sqrt{n}}\int_0^\vartheta\int_u^t M_{n,u}(s)(f_n(t,d(s-u),0)-g(t,d(s-u),0))dR(u)\nonumber\\
& = &\frac{1}{\sqrt{n}} \int_0^t\left[\int_0^{(t-w)\wedge
\vartheta} M_{n,u}(w+u)dR(u)\right](f_n(t,dw,0)-g(t,dw,0))\notag\\
& = &\frac{1}{\sqrt{n}} \int_0^t\tilde{M}_{n}(w,(t-w)\wedge\vartheta)\cdot (f_n(t,dw,0)-g(t,dw,0)),
\end{eqnarray}
where the last equality follows from the definition of  $\tilde{M}_{n}(w,\vartheta)$ in \eqref{tildeM}.
Then, by Proposition \ref{prop2},
 there exist partitions
 $0=w_0<w_1<\cdots<w_{N_0}=\tau$ and
$0=u_{n,0}\leq u_{n,1}\leq\cdots\leq u_{n,n_0}=\tau$ such that for all large $n$,
$$P\left(\sup_{i, j;
{w\in [w_i,w_{i+1})},\atop  u\in [u_{n,j},u_{n,j+1})}
\frac{1}{\sqrt{n}}| \tilde{M}_{n}(w,u) -
\tilde{M}_{n}(w_i,u_{n,j})|<\epsilon\right)>1-\epsilon.$$
Then, similarly to the derivation of (\ref{20}),
 we  have that for all large $n$, the following
holds with probability bigger than $1-2\epsilon$:
\begin{eqnarray}\label{22}
\eqref{proof2}
&\leq&
\frac{1}{\sqrt{n}}\sum_{i=1}^{N_0}\sum_{j=1}^{n_0}\left|\tilde
M_{n}(w_{i},u_{n,j})\int_{w_{i-1}}^{w_i}(f_n(t,dw,0)-g(t,dw,0))\right|
+ 2\epsilon K\nonumber\\
&\leq & 3\epsilon K.
\end{eqnarray}
Therefore the second term is also negligible. 
\end{proof}

\subsection*{A2: Proof of Proposition \ref{prop1}}\label{propproof1}
For the proof of Proposition \ref{prop1}, we shall make use of certain martingale inequalities as given in the following lemma, which is due to \cite*{Lenglart}.

\begin{lemma}\label{lenglart}
Let $\{W(s), 0\leq s\leq \tau\}$ be a square integrable martingale process whose sample paths are right
continuous with left limits. Then, for any $q>1$, there exists a
constant $C_q$ depending only on q, such that
\begin{eqnarray}\label{inequality1}
&&E\left(\sup_{s\leq \tau}|W(s)|^q\right)\leq C_q\left(E[\langle W
\rangle(\tau)]^{q/2}+E(\sup_{s\leq \tau}|\bigtriangleup
W(s)|^q)\right),
\end{eqnarray}
where $\langle W\rangle(s)$ denotes the predictable variation process of
the martingale  $\{W(s)\}$ and $\Delta W(s)=W(s)-W(s-).$

Moreover, if $\sup_{s\leq \tau}|\bigtriangleup W(s)|\leq c$, then for any $a,b>0$
\begin{eqnarray*}\label{inequality2}
 P\left(\sup_{s\leq \tau}|W(s)|\geq a, \langle W \rangle(\tau)\leq
b\right) \leq 2\exp\left(-\frac{a^2}{2b}\psi(ac/b)\right),\nonumber
\end{eqnarray*}
where $\psi(x)=2x^{-2}\{(1+x)[\log(1+x)-1]+1\}$.
\end{lemma}

\begin{proof}[Proof of Proposition \ref{prop1}]
Choose positive numbers $p$, $q>1$ such that $pq/2-p-q>1$. Let
$u_0=0$ and define $u_{n,j}$ recursively by
$$u_{n,j+1}=\inf\{\vartheta: \vartheta>u_{n,j}, 2\tau{\tilde K}_\tau(R(\vartheta)-R(u_{n,j})) \geq \epsilon^pn\}\wedge(u_{n,j}+\epsilon^p)\wedge\tau,$$
where $\tilde K_\tau$ is a constant satisfying
\begin{equation*}
\int_A\int_I q_n(ds\,du)<{\tilde K}_\tau\int_A\int_I
ds\,dR(u),\quad \hbox{ for any } A,I\subset [0,\tau].
\end{equation*}
It is easy to see from the above partition that there are at most $O(\epsilon^{-p})$
many, say $n_0$, distinct points in $[0,\tau]$. From
Lemma \ref{lemma1}, $\{W_{n,t,\vartheta}$, ${\cal F}_{n,t}$, $t\geq 0 \}$ is a
martingale, and we know that $u_{n,j}$, $j=1,\cdots,n_0$,  are
$\{{\cal F}_{n,t}, 0<t\leq \tau\}$ predictable. Thus,
$\{\sup_{\vartheta\in
[u_{n,j},u_{n,j+1}]}|W_{n,t,\vartheta}-W_{n,t,u_{n,j}}|, {\cal
F}_{n,t},t\geq 0\}$ is a nonnegative submartingale. By the Morkov
inequality and Doob's maximal inequality \citep{Doob},
\begin{eqnarray*}
&  & P\left( \max_{0\leq j<n_0} \sup_{\vartheta\in [u_{n,j},u_{n,j+1}];\atop 0\leq t\leq \tau}|W_{n,t,\vartheta}-W_{n,t,u_{n,j}}|\geq \epsilon\right)   \nonumber \\
& \leq & \frac{1}{\epsilon^q}\sum_{j=0}^{n_0-1} E\left( \sup_{\vartheta\in [u_{n,j},u_{n,j+1}];\atop 0\leq t\leq \tau}|W_{n,t,\vartheta}-W_{n,t,u_{n,j}}|^q\right) \nonumber\\
& \leq & \frac{1}{\epsilon^q}\sum_{j=0}^{n_0-1}
\left(\frac{q}{q-1}\right)^q E\left( \sup_{\vartheta\in
[u_{n,j},u_{n,j+1}]}|W_{n,\tau,\vartheta}-W_{n,\tau,u_{n,j}}|^q\right).
\end{eqnarray*}
Since $\{W_{n,\tau,\vartheta},{\cal
F}_{n,\tau,\vartheta},\vartheta\geq 0\}$ is a martingale and
$$
\sup_{\vartheta\in [u_{n,j},u_{n,j+1}]}
\bigtriangleup|W_{n,\tau,\vartheta}-W_{n,\tau,u_{n,j}}|\leq
\frac{1+\tilde K_\tau\tau}{\sqrt{n}},$$
 it follows from (\ref{inequality1}) that
\begin{eqnarray*}
&  & \frac{1}{\epsilon^q}\sum_{j=0}^{n_0-1} \left(\frac{q}{q-1}\right)^q E\biggr( \sup_{\vartheta\in [u_{n,j},u_{n,j+1}]}|W_{n,\tau,\vartheta}-W_{n,\tau,u_{n,j}}|^q\biggr)  \nonumber \\
& \leq & \frac{1}{\epsilon^q}\sum_{j=0}^{n_0-1} \left(\frac{q}{q-1}\right)^q C_q \left(E[\langle W_{n,\tau,u_{n,j}+\cdot}\rangle(u_{n,j+1}-u_{n,j})]^{q/2}+\frac{(1+\tilde K_\tau\tau)^q}{n^{q/2}}\right)\\
& \leq & C_q^*(\epsilon)^{pq/2-p-q}\leq {\epsilon},\nonumber
\end{eqnarray*}
where $C_q^*$ is a constant depending only on $q$ and the last
inequality holds when $\epsilon$ is sufficiently enough. Hence the desired
result follows.
\end{proof}

\subsection*{A3: Proof of Proposition \ref{prop2}}\label{propproof2}
To prove Proposition \ref{prop2}, we need the following lemma; see Lemma 5 in \cite{GuLai}.
\begin{lemma}\label{keylemma}
Let $q>0$ and $r>1$. Let $\{W_n, n\geq 1\}$ be a sequence of random variables defined in the same probability space and let $\{g_n\}$
be a sequence of nonnegative integrable functions on a measure space $({\cal X}, {\cal B}, \mu)$. Suppose that for every fixed $x\in{\cal X}$,
$g(x)$ is nondecreasing in $n\leq N$ and that
\begin{equation*}
E|W_i-W_j|^q\leq \left(\int_{\cal X}[g_i(x)-g_j(x)]d\mu(x)\right)^r
\textit{ for all } 1\leq j\leq i\leq N.
\end{equation*}
Then there exists a universal constant $C_{q,r}$ depending only on $q$ and $r$ such that
\begin{equation*}
E\left(\sup_{n\leq N}|W_n-W_1|\right)^q\leq C_{q,r}\left(\int_{\cal
X}[g_N(x)-g_1(x)]d\mu(x)\right)^r.
\end{equation*}
\end{lemma}

\begin{proof}[Proof of Proposition \ref{prop2}]
Choose positive numbers $p,q>1$ such that $pq/2-p-2q>1$. Let $w_0=0$,
and define $w_j$ recursively by $w_{j+1}=j \epsilon^p/{\tilde
K}_\tau$, where $\tilde K_\tau$ is a constant satisfying
\begin{equation*}
\int_A\int_I q_n(ds\,du)<{\tilde K}_\tau\int_A\int_I
ds\,dR(u), \hbox{ for any } A,I\subset [0,\tau].
\end{equation*} Denote $N_0=\floor{{\tilde K}_\tau\tau/\epsilon^p}+1$, and
redefine $w_{N_0}=\tau$.

Let $w_{n,i} = i\sqrt{\epsilon}/(nM)$ and  ${\cal N}_w=\{w_{n,i}:
i=0,1,\cdots,\floor{\tau Mn/\sqrt{\epsilon}}+1\}$ for some constant $M$. Then
$$P\left(\int_0^\tau\int_{u+w_{n,i}}^{u+w_{n,i+1}} p_n(du\,ds)\geq 2\right)=O(n^2)\cdot \epsilon/(nM)^2
\leq \epsilon/2$$
when $M$ is large enough.
By the definition of $\tilde K_\tau$, for
$\tilde{W}_{n,w,\vartheta}=\tilde{M}_{n}(w,\vartheta)/\sqrt{n}$, we have that
\begin{eqnarray}\label{first}
&      & P\left(\sup_{i, w_{n,i}\leq w\leq
w_{n,i+1}}|\tilde{W}_{n,w,\tau}-\tilde{W}_{n,w_{n,i},\tau}|\geq
2n^{-1/2}+{\tilde K}_\tau n^{-1/2}\right)\nonumber\\
& \leq & P\left(\sup_i \int_0^\tau\int_{u+w_{n,i}}^{u+w_{n,i+1}} p_n(du\,ds)\geq 2\right)\nonumber\\
&      & +P\left(\sup_i \int_0^\tau\int_{u+w_{n,i}}^{u+w_{n,i+1}} q_n(du\,ds)\geq {\tilde K}_\tau \right)\nonumber\\
&\leq  & \epsilon/2. 
\end{eqnarray}

Therefore, to prove Proposition \ref{prop2}, by (\ref{first}) and
the martingale property for $\{\tilde{W}_{n,w,\vartheta}, {\cal
F}_{n,\tau,\vartheta}, 0<\vartheta\leq \tau\}$ along entry time, we
only need to show that for any $\epsilon>0$,
$$P\left( \max_{0\leq j<N_0} \sup_{0\leq \vartheta\leq \tau \atop w\in
[w_j,w_{j+1}]\cap {\cal
N}_w}|\tilde{W}_{n,w,\vartheta}-\tilde{W}_{n,w_j,\vartheta}|\geq
\epsilon\right)<\epsilon/2,$$
for all large $n$.
Then, by Doob's inequality
and (\ref{inequality1}), similarly as in the proof of Proposition
\ref{prop1},
\begin{eqnarray*}
&      & P\left( \max_{0\leq j<N_0}  \sup_{0\leq \vartheta\leq \tau \atop w\in [w_j,w_{j+1}]\cap {\cal N}_w}|\tilde{W}_{n,w,\vartheta}-\tilde{W}_{n,w_j,\vartheta}|\geq \epsilon\right)   \nonumber \\
& \leq & \frac{1}{\epsilon^q}\sum_{j=0}^{N_0-1} E\left( \sup_{0\leq \vartheta\leq \tau \atop w\in [w_j,w_{j+1}]\cap {\cal N}_w}|\tilde{W}_{n,w,\vartheta}-\tilde{W}_{n,w_{j},\vartheta}|^q\right) \nonumber\\
& \leq & \frac{1}{\epsilon^q}\sum_{j=0}^{N_0-1}
\left(\frac{q}{q-1}\right)^q E\left( \sup_{w\in [w_j,w_{j+1}]\cap
{\cal N}_w}|\tilde{W}_{n,w,\tau}-\tilde{W}_{n,w_{j},\tau}|^q\right).
\end{eqnarray*}

For any $w_{n,i},w_{n,k} \in [w_j, w_{j+1}]\cap\, {\cal N}_w$, since
$\tilde{W}_{n,w_{n,k},\vartheta}-\tilde{W}_{n,w_{n,i},\vartheta}$
is a $\{{\cal F}_{n,\tau,\vartheta},\vartheta\geq 0\}$ martingale,
from (\ref{inequality1}) we have
\begin{eqnarray*}
 &     & E\left( |\tilde{W}_{n,w_{n,k},\tau}-\tilde{W}_{n,w_{n,i},\tau}|^q\right)\\
 &\leq & C_q \left(E\left[\langle \tilde{W}_{n,w_{n,k},\cdot}-\tilde{W}_{n,w_{n,i},\cdot}\rangle(\tau)\right]^{q/2}
 +\left(\frac{1+{\tilde K}_\tau(w_{n,k}-w_{n,i})}{n^{1/2}}\right)^q\right)\nonumber\\
 &\leq & C\epsilon^{-q/4} \left(\int_0^\tau \left[{\tilde K}_\tau\cdot 1(x\leq w_{n,k})-{\tilde K}_\tau\cdot 1(x\leq w_{n,i})\right]dx\right)^{q/2},\nonumber
\end{eqnarray*}
where $C$ is a constant. Then from Lemma \ref{keylemma},
there exists constant $C^*>0$ such that for all large $n$,
\begin{eqnarray}\label{last}
&     & \frac{1}{\epsilon^q}\sum_{j=0}^{N_0-1}
\left(\frac{q}{q-1}\right)^q E\left( \sup_{w\in
[w_j,w_{j+1}]\cap {\cal N}_w}|\tilde{W}_{n,w,\tau}-\tilde{W}_{n,w_{j},\tau}|^q\right)\nonumber\\
&\leq & \frac{1}{\epsilon^q}\sum_{j=0}^{N_0-1}\left(\frac{q}{q-1}\right)^q C\epsilon^{-q/4} \left(\int_0^\tau{\tilde K}_\tau\cdot 1(w_{n,i_{w_j}}<x\leq w_{n,i_{w_{j+1}+1}})\, dx\right)^{q/2}\nonumber\\
&\leq &  C^*(2\epsilon)^{pq/2-p-5q/4},
\end{eqnarray}
where $i_{w_j}=\max\{i: w_{n,i}\leq w_j\}$. By choosing $\epsilon$ sufficiently small, we have that the last term in (\ref{last}) must be smaller than
$\epsilon$. Hence the desired conclusion follows.
\end{proof}

\subsection*{A4: Lemma \ref{bgy}}
Lemma \ref{bgy} is used in the proof of Theorem \ref{consistency}. It is a restatement of Lemma A.5 in \cite{BGY}.

\begin{lemma}\label{bgy}
Consider a set of functions $\{f_{n,\alpha}: n\geq 1,\alpha\in A\}$
from $R^d$ to $R^d$. Suppose that (i) $\frac{\partial}{\partial
\theta}f_{n,\alpha}(\theta)$ are nonnegative definite for all $n$,
$\alpha$, $\theta$; (ii)
$\sup_\alpha|f_{n,\alpha}(\theta_0)|\rightarrow 0$ as
$n\rightarrow\infty$; (iii) there exists a neighborhood of
$\theta_0$, denoted by ${\cal N} (\theta_0)$, such that
$$\liminf_{n\rightarrow\infty}\inf_{\theta\in {\cal N}(\theta_0)} \inf_{a\in A} \lambda_{min}\left(\frac{\partial f_{n,\alpha}(\theta)}{\partial \theta}\right)>0,$$
where $\lambda_{min}$ is the minimum eigenvalue as defined in $C4$.
Then there exists $n_0$ such that for every $n>n_0$ and $\alpha\in
A$, $f_{n,\alpha}$ has a unique root $\theta_{n,\alpha}$ and
$\sup_{\alpha\in A}|\theta_{n,\alpha}-\theta_0|\rightarrow 0$.
\end{lemma}

\par

\noindent {\large\bf Acknowledgment}
The authors are grateful to the associate editor and two referees for their comments and suggestions, which lead to many clarifications and a better and more focused presentation. This research was supported in part by the NIH grant
5R37GM047845.
\par


\bibliographystyle{rss}
\bibliography{sequential}

\begin{thebibliography}{29}
\expandafter\ifx\csname natexlab\endcsname\relax\def\natexlab#1{#1}\fi
\expandafter\ifx\csname url\endcsname\relax
  \def\url#1{\texttt{#1}}\fi
\expandafter\ifx\csname urlprefix\endcsname\relax\def\urlprefix{URL }\fi

\bibitem[{Andersen \emph{et~al.}(1993)Andersen, Borgan, Gill and
  Keiding}]{ABGK}
Andersen, P.~K., Borgan, {\O}., Gill, R.~D. and Keiding, N. (1993).
  \emph{Statistical Models Based on Counting Processes}.
\newblock Springer, New York.

\bibitem[{Andersen and Gill(1982)}]{AndersonGill}
Andersen, P.~K. and Gill, R.~D. (1982). CoxÕs regression model for counting
  processes: a large sample study.
\newblock \emph{Ann. Statist.}, \textbf{10}, 1100--1120.

\bibitem[{{BHAT}(1982)}]{BHAT}
{BHAT} (1982). A randomized trial of propranolol in patients with acute
  myocardial infarction.
\newblock \emph{J. Amer. Med. Assoc.}, \textbf{147}, 1707--1714.

\bibitem[{Bilias \emph{et~al.}(1997)Bilias, Gu and Ying}]{BGY}
Bilias, Y., Gu, M. and Ying, Z. (1997). Towards a general asymptotic theory for
  cox model with staggered entry.
\newblock \emph{Ann. Statist.}, \textbf{25}, 662--682.

\bibitem[{Cox(1972)}]{Cox}
Cox, D.~R. (1972). Regression models and life-tables.
\newblock \emph{J. R. Statist. Soc. B}, \textbf{34}, 187--220.

\bibitem[{Cox(1975)}]{Cox1975}
Cox, D.~R. (1975). Partial likelihood.
\newblock \emph{Biometrika}, \textbf{62}, 269--276.

\bibitem[{Doob(1953)}]{Doob}
Doob, J.~L. (1953). \emph{Stochastic Processes}.
\newblock Wiley, New York.

\bibitem[{Fleming and Harrington(1991)}]{FlemingHarrington}
Fleming, T.~R. and Harrington, D. (1991). \emph{Counting Processes and Survival
  Analysis}.
\newblock Wiley, New York.

\bibitem[{Flournoy and Rosenberger(1995)}]{FR1995}
Flournoy, N. and Rosenberger, W.~F. (1995). \emph{Adaptive Designs}.
\newblock IMS, Hayward, CA.

\bibitem[{Gu and Lai(1991)}]{GuLai}
Gu, M.~G. and Lai, T.~L. (1991). Weak convergence of time-sequential censored
  rank statistics with applications to sequential testing in clinical trials.
\newblock \emph{Ann. Statist.}, \textbf{19}, 1403--1433.

\bibitem[{Hu and Rosenberger(2006)}]{HuRosenberger}
Hu, F. and Rosenberger, W.~F. (2006). \emph{The Theory of Response-adaptive
  Randomization in Clinical Trials}.
\newblock Wiley, New York.

\bibitem[{Kalbfleisch and Prentice(2002)}]{KalbfleischPrentice}
Kalbfleisch, J.~D. and Prentice, R.~L. (2002). \emph{The Statistical Analysis
  of Failure Time Data}.
\newblock Wiley, New York.

\bibitem[{{Kim, E. S., {\it et al.}}(2011)}]{Kim2011}
{Kim, E. S., {\it et al.}} (2011). The {BATTLE} trial: Personalizing therapy
  for lung cancer.
\newblock \emph{Cancer Discovery}, \textbf{1}, 44--53.

\bibitem[{Lan and DeMets(1983)}]{LanDeMets}
Lan, K. K.~G. and DeMets, D.~L. (1983). Discrete sequential boundaries for
  clinical trials.
\newblock \emph{Biometrika}, \textbf{70}, 659--663.

\bibitem[{Lenglart \emph{et~al.}(1980)Lenglart, Lepingle and
  Pratelli}]{Lenglart}
Lenglart, E., Lepingle, D. and Pratelli, M. (1980). \emph{Presentation unifiee
  de certaines in\'{e}gualit\'{e}s de la th\'{e}orie des martingales,
  S\'{e}minaire de Probabilit\'{e}s. Lecture Notes in Math,14,26-48}.
\newblock Springer, Berlin.

\bibitem[{Martinussen and Scheike(2000)}]{MartinussenScheike}
Martinussen, T. and Scheike, T.~H. (2000). A nonparametric dynamic additive
  regression model for longitudinal data.
\newblock \emph{Ann. Statist.}, \textbf{28}, 1000--1025.

\bibitem[{O'Brien and Fleming(1979)}]{O'BrienFleming}
O'Brien, P.~C. and Fleming, T.~R. (1979). A multiple testing procedure for
  clinical trials.
\newblock \emph{Biometrics}, \textbf{35}, 549--556.

\bibitem[{Pocock(1977)}]{Pocock}
Pocock, S.~J. (1977). Group sequential methods in the design and analysis of
  clinical trials.
\newblock \emph{Biometrika}, \textbf{64}, 191--199.

\bibitem[{Rebolledo(1980)}]{Rebolledo}
Rebolledo, R. (1980). Central limit theorem for local martingales.
\newblock \emph{Z.Wahr. verw. Geb.}, \textbf{51}, 269--286.

\bibitem[{Robbins(1952)}]{Robbins}
Robbins, H. (1952). Some aspects of the sequential design of experiments.
\newblock \emph{Bull. Amer. Math. Soc.}, \textbf{58}, 527--535.

\bibitem[{Sellke and Siegmund(1983)}]{SellkeSiegmund}
Sellke, T. and Siegmund, D. (1983). Sequential analysis of the proportional
  hazards model.
\newblock \emph{Biometrika}, \textbf{70}, 315--26.

\bibitem[{Siegmund(1985)}]{Siegmund1985}
Siegmund, D. (1985). \emph{Sequential Analysis: Tests and Confidence
  Intervals}.
\newblock Springer, New York.

\bibitem[{Slud(1984)}]{Slud}
Slud, E.~V. (1984). Sequential linear rank tests for two-sample censored
  survival data.
\newblock \emph{Ann. Statist.}, \textbf{12}, 551--571.

\bibitem[{Thompson(1933)}]{Thompson}
Thompson, W.~R. (1933). On the likelihood that one unknown probability exceeds
  another in view of the evidence of the two samples.
\newblock \emph{Biometrika}, \textbf{25}, 285--294.

\bibitem[{Wei(1978)}]{Wei}
Wei, L.~J. (1978). The adaptive biased coin design for sequential experiments.
\newblock \emph{Ann. Statist.}, \textbf{6}, 92--100.

\bibitem[{Wei and Durham(1978)}]{WeiDurham}
Wei, L.~J. and Durham, S. (1978). The randomized play-the-winner rule in
  medical trials.
\newblock \emph{J. Amer. Statist. Assoc.}, \textbf{73}, 840--843.

\bibitem[{Zelen(1969)}]{Zelen}
Zelen, M. (1969). Play the winner and the controlled clinical trial.
\newblock \emph{J. Amer. Statist. Assoc.}, \textbf{64}, 131--146.

\bibitem[{Zhang and Rosenberger(2007)}]{ZR2007}
Zhang, L. and Rosenberger, W.~F. (2007). Response-adaptive randomization for
  survival trials: the parametric approach.
\newblock \emph{Appl. Statist.}, \textbf{56}, 153--165.

\bibitem[{Zhou \emph{et~al.}(2008)Zhou, Liu, Kim, Herbst and Lee}]{Zhou2008}
Zhou, X., Liu, S., Kim, E.~S., Herbst, R.~S. and Lee, J.~J. (2008). Bayesian
  adaptive design for targeted therapy development in lung cancer - a step
  toward personalized medicine.
\newblock \emph{Clinical Trials}, \textbf{5}, 181--193.

\end{thebibliography}

\vskip .65cm
\noindent
Xiaolong Luo
\vskip 2pt
\noindent
Celgene Corporation 
\vskip 2pt
\noindent
E-mail: xluo@celgene.com
\vskip 4pt

\noindent
Gongjun Xu and Zhiliang Ying
\vskip 2pt
\noindent
Department of Statistics, Columbia University
\vskip 2pt
\noindent
E-mail: gongjun, zying@stat.columbia.edu
\vskip .3cm 



\end{document}